\documentclass[12pt]{article}
\usepackage[centertags]{amsmath}
\usepackage{amsfonts}
\usepackage{dsfont}
\usepackage{amssymb}
\usepackage{amsthm}
\usepackage{amsmath,mathrsfs}
\usepackage{amsmath,amscd}
\usepackage{mathtools}
\usepackage[matrix,arrow,curve,cmtip]{xy}
\title{\textbf{Consciousness in a Higher Categorical Context}}
\author{Renaud Gauthier \footnote{2020 Math. Subj. Class: 18N60, 81P10, 81P45. Keywords: Segal topos, micro-reversibility, information theory, Grothendieck construction, consciousness } \\ \\}
\theoremstyle{definition}
\newtheorem{Kay}{Remark}[subsection]
\newtheorem{inc}{Remark}[subsection]

\DeclareMathOperator*{\colim}{\text{colim}}

\newcommand{\beq}{\begin{equation}}
\newcommand{\eeq}{\end{equation}}

\newcommand{\rarr}{\rightarrow}

\newcommand{\eset}{\emptyset}
\newcommand{\Ob}{\text{Ob\,}}

\newcommand{\xrarr}{\xrightarrow}
\newcommand{\xlarr}{\xleftarrow}

\newcommand{\cC}{\mathcal{C}}

\newcommand{\cD}{\mathcal{D}}

\newcommand{\cG}{\mathcal{G}}

\newcommand{\cK}{\mathcal{K}}

\newcommand{\cT}{\mathcal{T}}

\newcommand{\cX}{\mathcal{X}}

\newcommand{\bK}{\mathbb{K}}

\newcommand{\bR}{\mathbb{R}}

\newcommand{\Cat}{\text{Cat}}

\newcommand{\Fun}{\text{Fun}}

\newcommand{\Hom}{\text{Hom}}

\newcommand{\Map}{\text{Map}}

\newcommand{\Set}{\text{Set}}

\newcommand{\uHom}{\underline{\Hom}}

\newcommand{\Catinf}{\Cat_{\infty}}

\newcommand{\dAff}{\text{dAff}}

\newcommand{\dAffCtildetau}{\dAff_{\cC}^{\;\sim  ,  \tau}}

\newcommand{\dSt}{\text{dSt}}

\newcommand{\eps}{\epsilon}

\newcommand{\oP}{\oplus}

\newcommand{\oT}{\otimes}

\newcommand{\RuHom}{\bR \uHom}

\newcommand{\skMod}{\text{s}k\text{-Mod}}

\newcommand{\SeT}{\text{SeT}}

\newcommand{\nuepsa}{\nu_{\eps}(a)}
\newcommand{\bOne}{\mathds{1}}
\newcommand{\dOne}{\text{d}\bOne}
\newcommand{\Bun}{\text{-Bun}}
\newcommand{\TBun}{\cT \text{-Bun}}
\newcommand{\TqBun}{\cT_q \text{-Bun}}
\newcommand{\TBunp}{\cT \text{-Bun}_p}
\newcommand{\TqBunp}{\cT_q \text{-Bun}_p}
\newcommand{\Td}{\cT_{\centerdot}}
\newcommand{\TdBund}{\cT_{\centerdot} \text{-Bun}_{\centerdot}}
\newcommand{\Grpd}{\text{Grpd}}
\newcommand{\peps}{ p \eps}
\newcommand{\pponeeps}{(p\text{+1})\eps}
\newcommand{\TUoneUq}{T_{U_1 \cdots U_q}}
\newcommand{\coCartFib}{\text{coCartFib}}

\begin{document}
\maketitle
\begin{abstract}
	We provide two representations of the Segal category $\cX$ modeling natural phenomena, the first one being based on the concept of micro-reversibility, producing a long sequence $\Sigma$ of categories as a resolution of $\cX$, the second one providing graded categories cofibered in groupoids over the categories of $\Sigma$, using the concept of consciousness as impetus. We show those two representations are dual to each other.
\end{abstract}

\section{Introduction}
In \cite{RG2} we considered dynamics in the Segal topos of derived stacks $\dSt(k)= \dAffCtildetau$, initially introduced in \cite{RG} as a means to model natural phenomena. This however is based only on one symmetric monoidal model category $\cC = \skMod$, for $k$ a commutative ring. As in \cite{RG4} one should consider all the base symmetric monoidal model categories $\cC_1, \cdots, \cC_N$ that come into play in describing natural phenomena, to give rise to $\prod_{1 \leq i \leq N} \dAff_{\cC_i}^{\; \sim, \, \tau_i} = \cX \in \SeT$, a Segal topos. Next suppose we have $P$ worlds to account for in order to have a comprehensive description of universal phenomena, each world with its very own Segal topos $\cX_j$ modeling local phenomena, $1 \leq j \leq P$. Then all phenomena would be accounted for by starting with $\cX = \prod_{1 \leq j \leq P} \cX_j \in \SeT$.\\

We would like now to relate such an object (or its dynamics $\delta^{\infty}\cX$ more generally, with $\delta^{\infty} \cX = \colim_n \delta^n \cX$, and $\delta \cX = \RuHom(\cX,\cX)$), with what we call the ground state of being $\Omega$, since natural phenomena must originate therein. Actually natural phenomena occur in $\cX$, dynamics thereof in $\delta^n \cX$, $n \geq 1$, including supernatural phenomena more generally, so $\colim \delta^n \cX = \delta^{\infty} \cX := \bOne$ gives a complete description of all phenomena all things considered, that is it is a presentation of $\Omega$. Hence it is $\delta^{\infty}\cX = \bOne$ we seek to relate to $\Omega$ Unfortunately, the latter is metaphysically defined as being devoid of any concept, and of points in particular, while $\bOne$ itself is a highly structured object. Thus presently we seek a ``smoother", point-independent formalism for the sake of presenting $\bOne$ from the perspective of $\Omega$, and this ushers the concept of fluid category, something we will define in this work. This concept is fundamentally built on the concept of micro-reversibility, according to which events that are really close are for all practical purposes indistinguishable, hence for instance the fluid nature of neighborhoods whose points are equivalent to each other by micro-reversibility.\\

To be more precise, the desired connection between $\bOne$ and $\Omega$ first involves a deconstruction in the sense of \cite{RG3}, an information-theoretic decomposition that is. Thus one first focuses on the deconstruction $dec \bOne$ of $\bOne$. Once this is done, we observe that in the infinite continuum $\Omega$, nothing is really isolated and a point independent formalism is therefore de rigueur to describe it. However $\bOne$ (and consequently $d\bOne = dec \bOne$) are objects built on a point-based formalism, which subsumes individual events can be isolated. To establish a relation between $\dOne$ and $\Omega$, one therefore needs to put $\dOne$ and $\Omega$ on a same footing, and by this we mean working with a point independent resolution of $\dOne$. This will prompt us to introducing the concept of fluid category to describe the appropriate resolutions of $\bOne$, from which it will become apparent that dynamics within $\dOne$ is given by a continuous flow of information. This fluidity we claim is the fundamental nature of phenomena.\\

A fluid category is a (gradual) resolution of a structured object such as $\dOne$ into something of the same nature as $\Omega$, and the reason for needing such an incremental process is that this pervasive unity within $\Omega$ translates locally as micro-reversibility at the level of $\bOne$, or $\dOne$ for that matter, hence turning $\dOne$ into a point-independent object must be done gradually. \\

\newpage

Practically, we will progressively cover $d\bOne$ by micro-reversibility-induced point-less neighborhoods that will ultimately yield a point-independent presentation of $d\bOne$. Letting $\eps$ measure micro-reversibility information theoretically, we will cover $\dOne$ with neighborhoods of size $\eps$ within which all points are equivalent by micro-reversibility, giving rise to a cover $\dOne[\eps]$ of $\dOne$. We will iterate this phenomenon by way of a thickening map $\delta$ which gives rise to a long sequence:
\beq
\dOne \xrarr{\delta} \cdots \xrarr{\delta} \dOne[\peps] \xrarr{\delta} \dOne[\pponeeps] \xrarr{\delta} \cdots \nonumber
\eeq
which consists in a progressive graded transition from $\dOne$, at which level one just deals with an ordinary category in a point-based formalism, to a colimit where point-less neighborhoods cover $\dOne$ in full, the whole sequence providing a point independent resolution of $\dOne$. It is in this limit that the statement $\dOne = \Omega$ is now well-defined. That long sequence from categories to this infinite continuum defines what we call a fluid category.\\

On the other hand we can focus on consciousness, arguing that $\Omega$ is pure consciousness, part of which is lost in the transition to the physical realm, which is operated by a reconstruction $rec:\Omega \rarr \bOne$. This process introduces obscurations, pure consciousness becomes tainted as a result, and we are left with consciousness as it is colloquially understood, an imperfect apprehension of reality. If the reader is uncomfortable with a seemingly imprecise notion of pure consciousness, we can equivalently proceed from our perspective; consciousness at our level is obviously tainted by obscurations, and is clarified once those obscurations are progressively removed. From that standpoint removing all obscurations results in a consciousness that we characterize as being pure. Now the problem we had earlier with $\bOne$ being point-based persists, so a statement of the form $rec \Omega = r\Omega = \bOne$ has to be interpreted in a point-independent formalism. For that purpose we recycle the work on micro-reversibility-based resolutions done on $\dOne$, now written in a reconstructed format for the sake of $\bOne$. Once this is done, we throw consciousness into the picture. One considers all possible consciousnesses at every point of focus $a$ of $\bOne \in \Cat$. Those are presented as smooth towers $T_a$, cone-shaped objects, based at $a \in \bOne$, whose cross-sections represent local consciousness relative to $a$, and ending in $\Omega$ in the limit. However because micro-reversibility is a natural phenomenon one has to take it into account. Thus it is not only over points $a$ of $\bOne \in \Cat$ that we have towers, it is also over the thick points of $\bOne[\peps]$ for $p \geq 1$ (neighborhoods of points equivalent to some core point $a$ via less than $p$ iterations of micro-reversibility), resulting in thickened towers. For each $p$, we have a functor $\chi_p: \bOne[\peps] \rarr \Grpd$ that to a thick point $U$ associates $\cT_U = \{T_U\}$, $T_U$ a thickened tower. Collectively those functors provide a prismatic decomposition of $\bOne$ into towers of evolving consciousness. Operating a Grothendieck construction on each of these functors, one obtains a graded category cofibered in groupoids over $\bOne$ presented as a fluid category. We argue both formalisms, the micro-reversibility-based one and the other based on consciousness, are dual to each other. As an immediate consequence of this duality, consciousness as we understand it appears as a particular filter on pure consciousness.

\section*{Relation to other works}
Most scientific works revolving around the concept of consciousness are written from a humanistic perspective, and are very much constructive, for the most part. As pointed out in \cite{MM}, things in a higher context are mirrored at a local level in the mind, which means that to access consciousness one can very well study neural networks, and from that perspective there are numerous publications that rely on tools from homotopy theory, motives, algebraic topology, etc... a (non-exhaustive of course) list of noteworthy publications being \cite{MnMi}, \cite{BB}, \cite{G}, \cite{M}, \cite{MM}, \cite{Mi}, \cite{T}, \cite{V}, \cite{CMnMi}, \cite{SKRKA}, \cite{KS}, \cite{KS2}. The point of making allusion to those works is to indicate that our present use of homotopy theory and other higher categorical ideas to tackle the concept of consciousness is nothing new. However our emphasis is really on the metaphysics of consciousness, and on understanding how consciousness as we understand it arises from pure consciousness. From this perspective we adopt a universal point of view and we focus on the realization of natural phenomena from a ground state of being, as opposed to just working with natural sciences as backdrop. We can also mention \cite{C} that alludes to the duality between the ground state of Being, natural phenomena, and a universe of all knowledge $\cK$ that we will briefly discuss.

\section*{Notations}
We denote by $\mu R$ the term micro-reversibility, by $\delta \Phi$ the $\delta$-formalism, $\bK \Phi$ the $\bK$-formalism. $\cK$ designates the universe of all knowledge, $\Omega$ the ground state of Being. $\bK$ stands for consciousness.

\section{$\delta$-formalism - Fluid categories}

\subsection{Backdrop}
In the introduction we gave two equations $r\Omega = \bOne$ and $d \bOne = \Omega$, that are to be understood in a point-independent formalism. Part of the work that follows will revolve around transforming $\bOne$ where natural phenomena take place, into its point-independent instantiation. Along those lines, recall that in  \cite{RG3} we introduced $IRV$s, gradual, geometric, complete information theoretic decompositions of concepts of interest, presumably objects of $\bOne$. $IRV$s are objects of a category that we denoted $\cG_0'$, with an isomorphic category $\cG_0$ whose objects are reversed $IRV$s seen as reconstructions of concepts from bits of information in $\Omega$. $IRV$s and their inverses being maps in themselves we have $\Ob \cG_0 \subset \Map(\Omega, \bOne)$, and $ \Ob \cG'_0 \subset \Map(\bOne,\Omega)$. At that point in \cite{RG3} we also mentioned that to form $IRV$s, one has to put oneself from the perspective of an all-encompassing Source that has a universe $\cK$ of all knowledge at its disposal for the sake of constructing $IRV$s, since those are complete, which presupposes having access to all information. Since the entirety of $\cK$ is used in constructing $IRV$s, the argument being that information in $\cK$ is necessarily used to produce objects, which have an $IRV$ of their own, it follows one can write:
\beq
\cG'_0 \doteq \cK \doteq \cG_0 \nonumber
\eeq
and having $\Ob \cG_0 \subset \Map(\Omega, \bOne)$ and $\Ob \cG'_0 \subset \Map(\bOne, \Omega)$ one can represent the above representation of $\cK$ diagrammatically as:
\beq
\xymatrix{
	&\cK \ar[dl] \ar[dr] \\
	dec: \bOne \rarr \Omega  && rec: \Omega \rarr \bOne 
} \nonumber
\eeq
where the $rec/dec$ maps are being understood $IRV$-style for instance. Compactly present this like so:
\beq
\xymatrix{
	& \cK  \\
	\Omega \ar@/^2pc/[rr] && \bOne \ar@/_2pc/[ll]
} \label{arch}
\eeq
This simple heuristic presentation shows $\Omega$ is not the Source, as one would naively think. The Source itself is represented as $\cK$ information-theoretically. Indeed it is defined to be that object from which everything that exists originates, it is characterized by existence, which is pure energy, i.e. it has dynamical properties (\cite{RG5}) as exemplified by the above diagram, if we tentatively accept that $\cK$ is a feature of the Source. To justify this fact energy is also information (\cite{J1}, \cite{J2}, \cite{Lt}, \cite{S1}, \cite{S2}), or equivalently knowledge, thus the Source itself can be represented as $\cK$. When in \cite{RG3} we refer to $\cG'_0$, or $\cG_0$, as being the ground state of natural phenomena, this is to be understood in a dynamical sense, case in point the fact that $\Ob \cG_0 \subset \Map(\Omega, \bOne)$ for instance. In that paper one showed we have structured objects and their decompositions coexisting simultaneously with $IRV$s bridging the gap between those representations, thus we can take $IRV$s themselves to provide a dual presentation of natural phenomena by their ground state decompositions, hence the name given to $\cG_0$ and $\cG_0'$. Here calling $\Omega$ the \textbf{ground state of Being} that's in contrast to having $\bOne$ on the other side of \eqref{arch}. Observe that the picture \eqref{arch} above shows we have simultaneity of $\Omega$ and $\bOne$, they coexist at the level of the source $\cK$.

\subsection{Micro-reversibility}
\subsubsection{micro-reversibility relation}
It is well-known that there is such a physical concept as \textbf{micro-reversibility} \cite{L}, \cite{B}, which we will abbreviate by $\mu R$. We denote by $\sim$ the relation characterizing $\mu R$. When we write $a \sim b$ we mean exactly this, that $a$ and $b$ are interchangeable. This relation is reflexive and symmetric. At the level of $\Omega$, it is transitive for lack of reference points. \\

Strictly speaking, physical micro-reversibility for an arrow $a \xrarr{f} b$ means $f$ has a left inverse $ a \xlarr{f^{-1}} b$. We claim the fundamental nature of $\mu R$ is not local reversibility of dynamic flows but stems from the underlying unification of phenomena at the information theoretic level, which we refer to as (abstract) $\mu R$.

\subsubsection{$\mu R$ and information theory}
When one applies $\mu R$, this is information-wise as we now argue. There is a relative notion of $\mu R$: $a$ and $b$ objects of a category $\cC$ are $\cC$-$\mu R$ if for all $c \in \cC$, $\Hom_{\cC}(c,a) \xrarr{f_*} \Hom_{\cC}(c,b)$ is an isomorphism for all $f:a \rarr b$. We say $a$ and $b$ are micro-reversible if they are $\cC$-micro-reversible for all categories $\cC$ they can be seen as being objects of. To put all those categories on a same footing, we must work at the information-theoretic level. In any case, identifying objects in this manner is a characteristic of $\Omega$, and it is from $\Omega$ that we initially get this notion.\\ 

From the perspective of a completely deconstructed environment such as $\Omega$ (in the information theoretic sense), one can clearly identify what are the events that are close to each other in the absolute, all dimensions/characteristics combined. Now as just said, all ``points" in $\Omega$ are transitively equivalent, it being an infinite continuum, thus strictly speaking one does not need a physical concept of micro-reversibility. When $\Omega$ is manifested in the natural realm $\bOne$ (or equivalently $\dOne$) with its persistent point-based character, this global equivalence of $\Omega$ is noticeable only locally on a very small scale, and this is what we refer to as micro-reversibility. Because $\mu R$ is a local manifestation of the equivalence of information theoretic objects in $\Omega$, we take $\mu R$ to be determined by information theoretic content.

\subsubsection{Construction of $\mu R$-neighborhoods}
That $\bOne$ and $\dOne$ manifest such a thing as $\mu R$ only, as opposed to a global equivalence as in $\Omega$, this means those spaces are considered to be point-based, the reason being that $\mu R$ is a local phenomenon, relative to a fixed point of interest $a$. Presently we focus on $\dOne$ in its point-based incarnation. $\dOne$ is $\bOne$ in raw deconstructed form, within an ambient vista $\Omega$ (\cite{RG3}). Regarding $\Omega$ as a relative vista will allow us to ultimately make the transition to the statement $\dOne = \Omega$, the reason being that locally within $\Omega$, objects of $\dOne$ will see what are neighboring events in accordance with $\mu R$, and this we repeat until the point $\dOne$ is presented in a completely point-independent formalism, at which point $\dOne = \Omega$ in full. We denote by $\eps$ the information-theoretic size that defines $\mu R$, something we can regard as a quantum cutoff. In other terms in the natural realm, two events $a$ and $b$ are regarded as indistinguishable by $\mu R$ if $d(a,b) < \eps$, where $d$ here measures informational content. We denote by $\nuepsa$ the neighborhood of a point $a \in \dOne$, within which all points are equivalent to $a$ by $\mu R$, that is all those points $b$ such that $d(a,b) < \eps$. We regard this as a \textbf{thickened point} $a$. For a core point $a$, points of $\nuepsa$ are therefore those concepts that are closest to $a$ information-theoretically, since $\mu R$ is an information theoretic concept.

\subsubsection{Miscellaneous remarks about $\mu R$-neighborhoods}
Since in the natural realm phenomena are dynamic, $\mu R$ is observed using time, hence points of $\nuepsa$ other than $a$ are pre-existing, they are yet to occur. In contrast to \cite{RG2} where q-deformations were defined using differential matrices, which didn't include higher morphisms, presently in a deconstructed setting $\dOne$ with $\bOne = \delta^{\infty}\cX$ all information is considered, higher morphisms included, hence those neighborhoods are all inclusive. Along those same lines, since $\Omega$ corresponds to pure consciousness as we will see, insofar as there are no obscurations therein, points of $\nuepsa$ correspond to regions of total awareness relative to their center $a$.\\

\subsection{$\mu R$-cover of $\dOne$}
\subsubsection{$\eps$-cover of $\dOne$}
In this subsection we argue one can cover the point-based object $\dOne$ by point-less neighborhoods of the form $\nuepsa$, which provides a local immersion of $\dOne$ into $\Omega$. An immediate upshot of having this cover is that one can then talk about points of $\Omega$ without contradicting the fact that there are no points in $\Omega$, by just invoking the relation $\dOne = \Omega$. Locally one can focus on points $a \in \dOne$ within the vista $\Omega$. \\

The global equivalence in $\Omega$ is reflected by an $\eps$-local $\mu R$ in $\dOne  $. Small neighborhoods of $\mu R$-equivalent points are denoted by $\nuepsa$ above. Such thick points of $\dOne$ provide a cover thereof. We can write:
\beq
\dOne \doteq \bigcup_{a \in \dOne} \nuepsa := \dOne[\eps] \nonumber
\eeq
where $\dOne[\eps]$ is the collection of all thick points of $\dOne$, in accordance with $\mu R$.

\subsubsection{$\dOne[\eps]$ category}
We turn $\dOne[\eps]$ into a category as follows. Focusing on individual thick points, first a given point $a \in \dOne$ is immersed in $\Omega$, and within $\Omega$ we have a perfect arrangement of deconstructed objects. This allows us to populate $\eps$-neighborhoods $\nuepsa$. Those are $\eps$-immersions of $\dOne$ into $\Omega$, they are worlds unto themselves. Thus from the perspective of a point $a$, the core of some element $U = \nuepsa$ of $\dOne[\eps]$, the only events that it can notice are those of $U$ by $\mu R$. For $V$ another thick point of $\dOne[\eps]$ distinct from $U$, it does not make sense for elements of $U$ to spontaneously travel to $V$. Thus any morphism of $\dOne[\eps]$ between distinct thick points must be defined relative to the infinite continuum $\Omega$ viewed as a vista. In $\Omega$ one can focus on one thick point or another, those are well-defined regions within $\Omega$. This change of focus defines morphisms between objects of $\dOne[\eps]$. One can represent this diagrammatically as follows:
\beq
\xymatrix{
	& \Omega \ar[dr]^{\pi_V} \\
	U \ar[ur]^{a_U} \ar[rr]_f && V
} \nonumber
\eeq
where $a_U$ will be formally defined later. For the moment we take this to be an inclusion of $U$ into $\Omega$ following the vista line of thought. $\pi_V$ amounts to taking a projection from $\Omega$ down to $V$. The definition is simple enough that we can see we have an obvious identity and composition of maps. Since we have $a_U \circ \pi_U = id_{\Omega}$, we have composition of morphisms that we represent diagrammatically by:
\beq
\xymatrix{
	& \Omega \ar[dr]_{\pi_V} \ar@{=}[rr] && \Omega \ar[dr]^{\pi_W} \\
	U \ar[ur]^{a_U} \ar[rr]_f && V \ar[ur]^{a_V} \ar[rr]_g && W
} \nonumber
\eeq
being compressed to:
\beq
\xymatrix{
	&\Omega \ar[dr]^{\pi_W} \\
	U \ar[ur]^{a_U} \ar[rr]_{g \circ f} && W
} \nonumber
\eeq
with the intermediate $\pi_V \circ a_V$ zigzag in between $\Omega$ and $V$ being understood. To put $\Omega$ at the top of this triangular diagram in spite of the fact that thick points are objects of $\dOne[\eps]$ is to emphasize that in the limit when neighborhoods grow in size as we will do next, we work with $\Omega = \dOne $ in a point-independent setting. Presently, $U$ and $V$ are just objects of $\dOne[\eps]$, with $\Omega$ as ambient vista. Nevertheless we have a slight improvement over simply working with $\dOne$, since in writing $\dOne \doteq \dOne[\eps]$ we can now instead work with $\dOne[\eps]$, which is still physically based at points, but with objects $\eps$-neighborhoods defined via $\mu R$, which is the nature of $\Omega$. This is a first step towards turning $\dOne$ into a point-independent object $\dOne = \Omega$.

\subsection{$\delta$-formalism}
\subsubsection{$\delta$ map}
If one can consider $\mu R$ at the level of single points $a \in  \dOne$ within the vista $\Omega$, points of $\nuepsa$ that are equivalent to $a$ can formally be considered for $\mu R$ as well. However at the level of natural phenomena $\mu R$ is not transitive. This means if $a \sim b$, and $b \sim c$, in such a manner that the deconstructed information theoretic contents of $a$ and $c$ differ by more than $\eps$, it is not true that $a \sim c$ since $\eps$ measures the cutoff for $\mu R$ in the natural realm. This means points equivalent to points of $\nuepsa - a$ by $\mu R$ that fall outside of $\nuepsa$ exist only as mere potentialities relative to $a$. One can formally illustrate this fact by introducing a grading on neighborhoods, whereby $\nuepsa = \nuepsa^{(0)}$ is a neighborhood of degree zero insofar as points therein are immediately accessible to $a$ by $\mu R$, whereas for $b \in \nuepsa$, $b \neq a$, some points of $(\nu_{\eps}(b) - \nuepsa)$ are not yet accessible to $a$ by $\mu R$, so are defined to be of degree one relative to $a$, and we would denote this by $(\nu_{\eps}(b) - \nuepsa)^{(1)}$. Thus:
\begin{align}
	\nu_{\eps}(b) &= (\nu_{\eps}(b)-\nu_{\eps}(a)) \cup (\nuepsa \cap \nu_{\eps}(b)) \nonumber \\
	&= (\nu_{\eps}(b)-\nu_{\eps}(a))^{(1)} \cup  (\nu_{\eps}(b)-\nu_{\eps}(a))^{(0)} \cup (\nuepsa \cap \nu_{\eps}(b)) \nonumber \\
	&= (\nu_{\eps}(b)-\nu_{\eps}(a))^{(1)} \cup ( \nuepsa \cap \nu_{\eps}(b)) \nonumber 
\end{align}
since $(\nu_{\eps}(b) - \nuepsa)^{(0)} = \eset$. Now formally applying $\mu R$ to all points of $\nuepsa$, we get a larger neighborhood of $a$, which becomes something we denote by
\begin{align}
	\nu_{2 \eps}(a) &= \bigcup_{b \in \nuepsa} \nu_{\eps}(b) \in \dOne[2 \eps] \nonumber \\
	&= \nuepsa \cup ( \bigcup_{b \in \nuepsa ,\;  b \neq a} \nu_{\eps}(b) - \nuepsa)^{(1)} \nonumber
\end{align}
Note that points of $\nu_{2 \eps}(a) - \nuepsa$ are of degree one relative to $a$. In the same manner that we saw $\dOne[\eps]$ as a category, we define $\dOne[2 \eps]$ with objects $\nu_{2\eps}(a)$, $a \in \dOne$, as a category. We will later define what its morphisms are. This operation of thickening $\eps$-neighborhoods we just performed we denote by $\delta$:
\beq
\delta: \dOne[\eps] \rarr \dOne[2 \eps] \nonumber
\eeq
This can obviously be repeated, and we obtain a long sequence:
\beq
\cdots \xrarr{\delta} \dOne[ \peps] \xrarr{\delta} \dOne[\pponeeps] \xrarr{\delta} \cdots \nonumber
\eeq
where $\nu_{\peps}(a) \in \dOne[\peps]$ can be decomposed as:
\beq
\nuepsa \cup (\nu_{2\eps}(a) - \nuepsa) \cup \cdots \cup (\nu_{\peps}(a) - \nu_{(p-1)\eps}(a)) \nonumber
\eeq
with $\nu_{k\eps}(a) - \nu_{(k-1)\eps}(a)$ of degree $k-1$.

\subsubsection{Remarks about the directed system}
Each stage in this directed system amounts to enlarging our horizon within the ambient vista $\Omega$, so that in the limit, which is formally a colimit, one recovers $\Omega$ in full as we will show, i.e. $\colim_n \delta^n \dOne = \Omega = \dOne$, a point-independent statement.\\

Observe that to hop from one category of thick points to the next via $\delta$ gives us a quantization of this directed system that is imposed on us by $\mu R$.

\subsubsection{$\peps$-neighborhoods interpretation}
Note that as before any $\nu_{\peps}(a)$ can be regarded as a thick point centered at $a$, with points therein further away from it to appear as faded to pictorially represent this grading of the different strata of $\nu_{\peps}(a)$ . We can picture this as a concentric information-theoretic neighborhood with a core point, with a degree-wise fading out as we move away from the core point $a$, having the equivalent interpretation that points further away from $a$ are less and less likely to occur, the reason being that they are further removed from $a$ information theoretically hence are more and more hypothetical. Points closer to the center are having a strong bond with $a$. This picture can also be regarded as providing a geometric representation of probabilistic events relative to $a$.\\

\subsubsection{Colimit of the directed system}
Let $U \in \dOne[\eps]$ centered at some point $a$, and denote by $\iota$ the map that gives us the information content of an object. Observe that in the limit when $U$ is really large as it progresses along the directed system, $\iota(a) \ll \iota(U)$. At that point one is effectively losing track of what core point $a$ we started with, that is there is no longer any referent, from which it follows $\mu R$ becomes transitive. Additionally $\colim \delta^n U$ covers $\Omega$. These two facts combined lead to $\colim \delta^n U = \Omega$. Doing this for all $U \in \dOne[\eps]$ it follows:
\beq
\colim \delta^n \dOne[\eps] = \Omega \label{colim}
\eeq
The advantage of having this picture of a directed system is that $\Omega$ is inert, with no manifest information theoretic dynamics, and no referents. However writing $\Omega = \colim \delta^n \dOne[\eps]$ provides us with a dual presentation of $\Omega$ where information flows relative to points can be defined. Indeed for a fixed point $a \in \dOne[\eps]$, the collection of neighborhoods $\nu_{\peps}(a)$ progressing along this system as in:
\beq
\cdots \xrarr{\delta} \nu_{p \eps}(a) \xrarr{\delta} \nu_{(p+1)\eps}(a) \xrarr{\delta} \cdots \nonumber
\eeq
with colimit $\Omega$ dynamically represents a gradual flow of information away from $a$, since in each $\nu_{p \eps}(a)$ points closer to $a$ have a stronger link with $a$, and points further away have a weaker connection.

\subsubsection{Fluid category definition}
We can actually have that long sequence start at $\dOne[0] = \dOne$, which is a category, thereby showing the long sequence:
\beq
\dOne \xrarr{\delta} \cdots \xrarr{\delta} \dOne[\text{p}\eps] \xrarr{\delta} \dOne[\text{(p+1)} \eps] \xrarr{\delta} \cdots \nonumber
\eeq
provides a $\mu R$-resolution of a category, namely $\dOne$. Because we start with $\dOne$, which is a (Segal) category, and the fact that points of $\nu_{\peps}(a)$ are equivalent by $\mu R$, so that we can regard those neighborhoods as being fluid neighborhoods, the entire resolution has for aim to reveal the fluid character of $\dOne$, complete in the limit, so the colimit of that sequence we call a \textbf{fluid category}. Alternatively, this resolution effectively gives a gradual local immersion of $\dOne$ into the ambient vista $\Omega$, so we identify this resolution with a dynamic presentation of $\dOne$ itself. Initially $\dOne \in \Cat$, in the colimit $\colim \delta^n \dOne$ is of the nature of $\Omega$, and the sequence in between provides us with a transition from a point-based formalism to a fully point-independent formalism. Contrast this with the work of \cite{RG2} where dynamics in $\cX$ took place in $\delta \cX = \RuHom(\cX,\cX)$, with objects $\psi$ therein providing one step only in the dynamics of $\cX$, with $\colim \psi$ providing a flow. Here in contrast we deal with fluid neighborhoods with no manifest discretization of the dynamics since it is built in $\bOne = \delta^{\infty}\cX$. In the limit we have effectively replaced $\dOne$ by its point-independent (fluid category) representation. But in the limit we also recover $\Omega$ by \eqref{colim}. Thus we have $\Omega = \dOne$ in a point-independent formalism. The formalism associated with this picture we refer to as the \textbf{$\delta$-formalism}, which we abbreviate by $\delta \Phi$.\\

\subsubsection{categorical formalism}
- \underline{Categories $\dOne[\peps]$} \\
We now go back to that claim above that each $\dOne[\peps]$ is a category. Having $\colim \delta^n U = \Omega$, denote taking the colimit of iterations of $U \in \dOne$ along the directed system by $a_U$, $a$ being short for ascendancy. This really means that we enlarge $U$ within the vista $\Omega$ till we reach the point we have become $\Omega$ properly speaking. The reverse process, focusing on $U$ from the perspective of $\Omega$ we denote by $\pi_U$. For simplicity we first start with the case $p=1$. $f:U \rarr V$ in $\dOne$ induces a map in $\dOne[2 \eps]$ defined by:
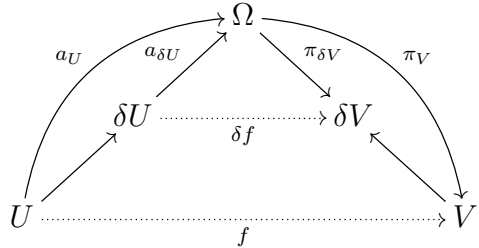
\beq
\xymatrix{
	&& \Omega \ar[dr]^{\pi_{\delta V}} \ar@/^2pc/[ddrr]^{\pi_V} \\
	& \delta U \ar[ur]^{a_{\delta U}} \ar@{.>}[rr]_{\delta f} && \delta V \\
	U \ar@/^2pc/[uurr]^{a_U} \ar[ur] \ar@{.>}[rrrr]_{f} &&&& V \ar[ul]
} \nonumber
\eeq
This has to be understood as follows: one can make sense of $\delta U$ and $\delta V$ individually by $\mu R$, by constructing them from $U$ and $V$ respectively using point-wise applications of $\delta$. Those are objects of $\dOne[2\eps]$, which originate from $\dOne[\eps]$ in the directed system. Now in a deconstructed environment, $U$, $V$ and $f$ are put on a same footing, so in the same manner that one can thicken $U$, $V$, etc..., if we have a map $f: U \rarr V$, by $\mu R$ one can thicken $f$ to $\delta f$ as well, in such a manner that it fits in the above diagram. Precisely, $\delta(U \cup f \cup V)$ is a thick point. Because $f$ has components other than those of $U$ and $V$, this thick point is made up of points that make up $\delta f$. Note that this means $\delta f$ is not a single morphism but is a \textbf{thickened morphism}. Morphisms of $\dOne[\peps]$ are constructed in this fashion more generally. Composition is defined as in $\dOne[\eps]$.\\

Observe that in such a picture, $f = \pi_V \circ a_U$, yet many such maps between $U$ and $V$ satisfy such a relation. This has to be understood as a limiting statement. Indeed, once $U$ (or $V$ for that matter) has been enlarged to the point it is the size of $\Omega$ information-theoretically, at that point there is no longer any map from $\colim \delta^n U$ to $\colim \delta^n V$, for the simple reason that both are $\Omega$, and we are only left with an identity between them.\\

-\underline{$\delta$ functor}\\
Coming back to the diagram above, by the above argument again, write $\delta f = \pi_{\delta V} \circ a_{\delta U}$ in the limit. We now argue $\delta$ is a functor on the categories $\dOne[\peps]$. For one thing we have $\delta id_U = id_{\delta U}$. To understand why that is the case, one simply uses the very definition of $\delta$; when this map is operating on an object $x$, it selects neighboring objects information-wise that are closest to $x$, meaning those objects that have the least amount of difference in information from that of $x$. Thus starting from some thick point $U \in \dOne[\peps]$, among all the possible maps that $\delta id_U$ could be, $id_{\delta U}$ is the one with the least amount of deviation from $id_U$ on the original object $U$, thus minimally we have $\delta id = id$. Additionally, thickened maps such as $\delta f$ must fit in a triangular diagram such as the one above, and the composition of such diagrams yields $\delta(g \circ f) = \delta g \circ \delta f$ by minimality, which we can represent as follows:
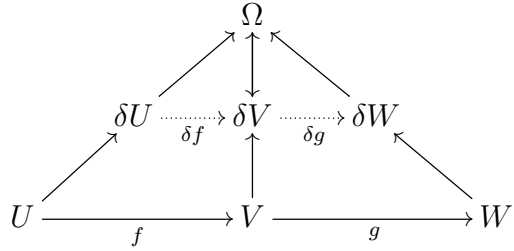
\beq
\xymatrix{
	&& \Omega \ar[d] \\
	& \delta U \ar[ur] \ar@{.>}[r]_{\delta f} & \delta V \ar[u] \ar@{.>}[r]_{\delta g} &\delta W \ar[ul] \\
	U \ar[ur] \ar[rr]_f && V \ar[u] \ar[rr]_g && W \ar[ul]
} \nonumber
\eeq
This is true at any level of the long sequence, showing that $\delta: \dOne[\peps] \rarr \dOne[\pponeeps]$ is a functor more generally.\\

\subsection{Natural $\delta$-formalism}
\subsubsection{Long sequence for $\bOne$}
$\dOne$ is $\bOne$ in deconstructed form, which is useful for the sake of constructing $\mu R$-neighborhoods within $\Omega$. However it is $\bOne$ that we observe in nature. The reconstructed version of the long sequence for $\dOne$ gives us:
\beq
\bOne \xrarr{\delta} \cdots \xrarr{\delta} \bOne[\text{p}\eps] \xrarr{\delta} \bOne[\text{(p+1)} \eps] \xrarr{\delta} \cdots \nonumber
\eeq
with $\delta$ having the same effect as on the deconstructed categories $\dOne[\peps]$. Observe however, and this is important, that this is sequentially merely a reconstructed version of the $\delta \Phi$ for $\dOne$:
\beq
\xymatrix{
	\dOne \ar[d]^{rec} \ar[r]^{\delta} & \cdots \ar[r]^{\delta} & \dOne[\peps] \ar[d]^{rec} \ar[r]^{\delta} & \cdots \\
	\bOne \ar[r]^{\delta} & \cdots \ar[r]^{\delta} & \bOne[\peps] \ar[r]^{\delta} & \cdots \\
} \nonumber
\eeq
where a neighborhood $U$ of $\bOne[\peps]$ is of the form $rec (\nu_{\peps}(a))$ for some point $a$. Accordingly, we refer to this reconstructed version as the natural $\delta \Phi$.\\

\subsubsection{Homotopy wave function within the natural $\delta \Phi$}
Thick points $U \in \bOne[\text{p}\eps]$ being reconstructed thick points of $\dOne[\text{p}\eps]$, they are still centered at some point $a \in \Omega$, but they are no longer concentric, they are more like blobs of varying intensity. Now if thick points in the deconstructed picture offered an obvious geometric picture of probabilistic events, thick points in the reconstructed picture being blobs, now one needs probability distributions to describe them in full. Observe that it is in this reconstructed $\delta \Phi$ that the homotopy wave function $\psi_{AG}$ introduced in \cite{RG2} is well-defined. There are several reasons for this. For one, the $\delta \Phi$ provides a pristine setting, without obscurations, where events can take place and can be studied purely probabilistically. Since $\psi_{AG}$ is defined probabilistically and describes phenomena, it is well-defined within the natural $\delta \Phi$. Since also $\psi_{AG}$ describes dynamics, read flow of information, and this is given in full by the natural $\delta \Phi$, this is yet another reason for seeing the latter as the canvas on which $\psi_{AG}$ is defined. Another way to see that the natural $\delta \Phi$ is the proper setting for discussing $\psi_{AG}$ is that the latter is an integral curve. For this it needs an initial condition and a directional vector. Thick points $U = \nuepsa$ do provide a point $a$, and neighboring points in accordance to $\mu R$ represent most probable next events as well, and $\psi_{AG}$ would correspond to a cross-section of such neighborhoods.

\subsubsection{Reconstituted fluid category $r \bOne$}
Speaking of dynamics, one can write:
\beq
rec \; \delta \Phi := \{ \cdots \xrarr{\delta} \bOne[\peps] \xrarr{\delta} \bOne[\pponeeps] \xrarr{\delta} \cdots \} = \bOne = r\Omega \nonumber
\eeq
since the natural $\delta \Phi$ achieves the transformation from the physical system $\bOne$ into its point independent version $\bOne = r \Omega$. This sequence provides a dynamic, fluid representation of $\bOne$. However, the original $\delta \Phi$ is the one that makes most sense since it is then that neighborhoods assemble naturally within $\Omega$, it provides a progressive immersion of $\dOne$ into the vista $\Omega$. The natural $\delta \Phi$ so far is important for defining the homotopy wave function $\psi_{AG}$. Later we will see that it is also important for the formalism of consciousness.

\section{$\bK$-formalism - model of consciousness}
Pure knowledge in $\Omega$ is reconstructed into various components in the natural realm, some being regarded as pertaining to the realm of consciousness, which for brevity we denote by $\bK$, and the rest being considered as potential obscurations insofar as they can affect $\bK$. This dichotomy is therefore purely human-centric. Note that at the moment we will tentatively take $\bK$ to mean what it colloquially means.\\

For $a \in \Ob \bOne$ a system with consciousness, $\bOne_{a/} = \{\bOne_{a/} \xrarr{\delta} \bOne_{a/}[\eps] \xrarr{\delta} \cdots \}$ in the ambient vista $r\Omega$ provides us with complete knowledge from the perspective of $a$, which will now be prismatically decomposed into its $\bK$ components. To formalize this we go as follows. We assume the consciousness of $a$ grows smoothly, until it ideally reaches $\Omega$ in the limit. We represent different possible paths of consciousness starting at $a$ as towers over that point whose cross sections are open sets with areas of varying intensity, much as in $\bOne[p\eps]$, the difference being that thick points of $\bOne[\peps]$ contain pure knowledge whereas cross sections $U$ of towers are regions of incomplete knowledge since they correspond to a local $\bK$ relative a fixed point $a$. Those towers grow in size from simply $a$, and each cross-section as we move up the towers represent the evolution of $a$'s consciousness. Those towers are cone-shaped, based at $a$, at the apex of which there is no consciousness but $a$ itself. We mention towers, plural, since there are many different possible smooth evolutions of $\bK$ relative to a given point, hence the necessity to consider different towers. One can regard those towers as inverse systems of local $\bK$s projecting down to a point, with $\Omega$ as a limit of each tower.\\

This however is just based at a point. However because of $\mu R$, which is a physical phenomenon, one should really work over thick points $U = \nuepsa \in \bOne[\eps]$, something we regard as a mathematical version of the wave-particle duality. Thus next we work with $\bOne[\eps]$ and we consider thick points $U = \nuepsa$ therein. Towers over thick points $U$ are denoted $T_U$, and will be thickened towers in a sense to be precised, with $\lim T_U = \Omega$. The towers being thickened means that above $a$ one still has all possible towers of its evolving $\bK$, and each such tower forms the core of a thickened tower. Those are formed as follows. If $U = \nuepsa$ is obtained from $a$ by $\mu R$, cross-sections of towers over $a$ must be $\mu R$-thickened as well. We proceed much as in the definition of thick points in $\dOne[\eps]$, where $\mu R$ is applied point-wise. For a cross-section of a tower $T_a$ at a point $a$, $\mu R$ is applied point-wise to all of its points. Stratifying those cross-sections one gets the thickened tower picture. Above a thick point $U = \nuepsa \in \bOne[\eps]$, one gets a thickened tower $T_U$. The collection of all such towers $T_U$ is denoted $\cT_U$, and forms a groupoid, where a morphism from one tower $T_U$ to another $T'_U$ amounts to a smooth reshaping. We do this for any point $a$ of $\bOne$. To emphasize that those groupoids form a complete system of consciousnesses at thick points, we can denote this by $\bK(U) = \cT_U$. Working over $\bOne[\eps]$ is natural by $\mu R$. One ought to work with all categories $\bOne[\peps]$ for completeness' sake, hence towers should really be defined over thick points of $\bOne[\peps]$ as well. This allows us to consider dynamics in the $\mu R$-direction by utilizing towers over thick points of $\bOne[\peps]$ for $p > 1$, where those towers are constructed by repeating the cross-sectional thickenings of towers we have performed above, resulting in thickened towers over $\bOne[\peps]$ that not only provide the smooth $\bK$-evolution of a core point $a$, but its $\mu R$-dynamics as well, something we argued is necessary for a probabilistic treatment of $\psi_{AG}$. Alternatively, relative to $\mu R$, towers present truncations of pure $\mu R$-neighborhoods of $a$ with respect to consciousness, a filter of some sort.\\

Next we will see this can be formalized by way of an induced functor from the $\delta \Phi$ to the $\bK$-formalism, which we regard as a prismatic decomposition of pure consciousness in $\bOne \subset r\Omega$ into $\bK$ represented by pointed towers as defined above. We will see later that there is a way to get the $\delta \Phi$ back from the $\bK$-formalism. \\

We will also show that towers can merge, which amounts to having overlapping consciousnesses. This will lead to considering a functor:
\beq
\chi_{p,q}:\bOne[p\eps]^q \rarr \Grpd \nonumber 
\eeq
which to any $(U_1,\cdots,U_q)$ associates the merged towers $T_{U_1} \otimes \cdots \otimes T_{U_q}$ in a sense to be precised. Doing a Grothendieck construction on each of those functors for all $p$ and $q$ we obtain categories cofibered in groupoids, which fit into a long sequence that corresponds to a complete model of $\bK$, and which we refer to as the \textbf{$\bK$-formalism}.

\subsection{Grothendieck construction}
From \cite{Lu}, recall how the Grothendieck construction works: we start with a functor $\chi: \cC \rarr \Grpd$ for a category $\cC$. One constructs a category of elements $\cC_{\chi}$ whose objects are of the form $(c,\xi)$ for $c \in \Ob \cC$, and $\xi \in \chi(c)$, and morphisms are of the form $(f,\alpha)$, with $f:c \rarr c'$ a morphism of $\cC$, and $\alpha: \xi \rarr \psi$ with $\xi \in \chi(c)$, $\psi \in \chi(c')$, defined as an isomorphism $\chi(f)(\xi) \rarr \psi$, which is actually short for:
\beq
\xymatrix{
	\xi \ar@{.>}[dr]_{\alpha} \ar[r]^-{\chi(f)} & \chi(f)(\xi) \ar[d]^{\cong} \\ 
	& \psi
} \nonumber
\eeq

\subsection{Construction of the category $\TBun$}

We consider the functor $\chi: \bOne[\eps] \rarr \Grpd$ which to thick points $U$ of $\bOne[\eps] = \cC$ associates the collection of towers $\cT_U \in \Grpd$, with $\cT_U = \{ T_U \}$, $T_U$ tower (with $\lim T_U = \Omega$) regarded as a $\bK$-resolved neighborhood. Thus $\chi$ gives all the possible ways to reach $\Omega$ with $U$ as starting point, by way of enlarging one's consciousness. That $\cT_U$ is a groupoid is by virtue of the fact that one can always reshape one tower $T_U$ into another one $T'_U$ over the same base $U$. $\cT_U$ being a groupoid, for $f:U \rarr V$ in the base, the map $\chi f: \cT_U \rarr \cT_V$ is a functor, which sends objects of $\cT_U$ to objects of $\cT_V$, that is objectwise it is a morphism of towers, defined cross-section wise the way morphisms were defined in the $\delta$-formalism, and is originally based at $f$. To be more precise, $f: U \rarr V$ map of thick points is well-defined within the ambient vista $\Omega$. Say $\chi f: T_U \rarr T_V$ over $f$ to fix notations. Each cross-section $\tilde{U}$ of the thickened tower $T_U$ is made up of points, each of which is part of a $\mu R$-neighborhood by thickening, which maps to a $\mu R$-neighborhood of a point that is part of a cross-section $\tilde{V}$ of what will be $T_V$. The manner in which this is performed is encoded in $\chi f$, generalization of $f$ at the level of thick points to the level of thick neighborhoods. One can represent this diagrammatically as:
\beq
\xymatrix{
	& \Omega  \\
	U \ar[ur]^{T_U} \ar[rr]_f && V \ar[ul]_{T_V} 
} \nonumber
\eeq
the entire triangle representing the map $\chi f:T_U \rarr T_V$. Note this is just how the functor $\chi f$ acts on $T_U \in \cT_U$. For $T_U \simeq T'_U \xrarr{\chi f} T'_V \simeq T_V$, we have a commutative diagram:
\beq
\xymatrix{
	U \ar@{=}[dd] \ar[dr]^{T_U}  \ar[rr]^f && V \ar[dl]_{T_V} \ar@{=}[dd] \\
	& \Omega \\
	U \ar[ur]_{T'_U} \ar[rr]_f && V \ar[ul]^{T'_V}
} \nonumber
\eeq
whose vertical faces are given by reshapings $\phi: T_U \simeq T'_U$ and $\cX f(\phi):T_V \simeq T'_V$. For composites:
\beq
\xymatrix{
	&\Omega \\
	U \ar[ur]^{T_U} \ar[r]_f & V \ar[u]_{T_V} \ar[r]_g &W \ar[ul]_{T_W}
} \nonumber
\eeq
which provides:
\beq
\xymatrix{
	&\Omega \\
	U \ar[ur]^{T_U} \ar[rr]_{g \circ f} && W \ar[ul]_{T_W}
}\nonumber
\eeq
We denote $\cC_{\chi}$ constructed following the Grothendieck construction on $\chi: \bOne[\eps] \rarr \Grpd$ by $\TBun$, it is cofibered in groupoids over $\bOne[\eps]$ with $\TBun \xrarr{\pi} \bOne[\eps]$ forgetful functor. Its objects are of the form $(U,T_U)$, with $U \in \bOne[\eps]$, $T_U \in \cT_U$ so that $\pi(U,T_U) = U$.

\subsection{Generalization: $\TBunp$}
We naturally have another type of functor $\chi_p: \bOne[\peps] \rarr \Grpd$, $U \mapsto \cT_U = \{ T_U \}$, $T_U$ $p$-thickened tower over $U \in \bOne[\peps]$. The Grothendieck construction on this functor leads to a category $\TBunp$ over $\bOne[\peps]$, cofibered in groupoids.

\subsection{$\Td$ as symmetric monoidal category}
Once one has towers, which correspond to dynamic, evolving consciousnesses, we consider their merging together. The reason for even considering this possibility is fairly evident; towers grow in size, and all end in $\Omega$, so clearly there is a point when they meet. Towers merging together form multi-towers.\\

For $U,U' \in \bOne[\text{p} \eps]$, $T_U \in \cT_U$ and $T_{U'} \in \cT_{U'}$, we define their merging together by $T_U \oT T_{U'} = T_U \bigcup T_{U'}$. Generalize this as follows:
\begin{align}
	\cT_U \oT \cT_{U'} &= \{ T_U \} \oT \{ T_{U'} \} \nonumber \\
	&= \{ T_U \oT T_{U'} \} \nonumber \\
	&= \{ T_U \cup T_{U'} \} \nonumber \\
	&= \cT_U \cup \cT_{U'} \nonumber
\end{align}
Now observe that towers provide dynamic consciousnesses; if $a$ is the core of $U = \nu_{\peps}(a)$, then a given tower $T_U$ over it gives us a path of increasing $\bK$ starting from $U$, thus $\cT_U = \{ T_U \}$ gives us all the possible perceptions from $a$ within $\bOne[\peps]$, which is reminiscent of the perception $\Hom(a,-)$. Now $\Hom(a,-)$ is just a $\Set$-based functor, while $\bK(a) = \cT_a$ is a geometric object that enhances $\Hom(a,-)$, with $\cT_a = \{ T_a \}$, $T_a$ a generalized $\bK$-neighborhood. We have $\cT_a \oT \cT_b \doteq \bK(a \oP b)$ since this product provides trees with two feet $a$ and $b$, which we represent as $a \oP b$. 
\beq
\bK(a \oP b) = \cT_a \oT \cT_b =  \bK(a) \oT \bK(b) = \bK(a) \cup \bK(b) \nonumber
\eeq
which generalizes:
\beq
\Hom(a \oP b, -) = \Hom(a,-) \oP \Hom(b,-) \nonumber
\eeq
To summarize, towers merging together consists in adding the perceptions of their core elements.\\

Having defined the product on towers, we are led to introducing $\Td = \oP_{n \geq 1} \cT_n$, with $\cT_n = \{T_{U_1 \cdots U_n}\} = \cup_{U_1 \cdots U_n} \cT_{U_1 \cdots U_n}$ where $T_{U_1 \cdots U_n} = T_{U_1} \oT \cdots \oT T_{U_n}$. Here it is understood that the feet $U_1, \cdots, U_n$ of a given tower are defined over a same base $\bOne[\peps]$. Formally we have $\oT: \cT_n \times \cT_q \rarr \cT_{\leq n+q} = \oP_{i \leq n+q} \cT_i$. The unit for the tensor product is the empty tower, over the empty thick point. By definition of $\oT$, symmetry and associativity are obvious, and we have the triangle diagrams for left and right multiplication. This makes $(\Td,\oT)$ into a symmetric monoidal category. To come back to the product $\oT$, even if it is true that the $\cT_n$ are constructed iteratively, given a generic element $T_{U_1 \cdots U_n}$ of $\cT_n$, one needs to make sense of its representation as $T_{U_1} \oT \cdots \oT T_{U_n}$. Observe that many towers fit the bill, namely all those towers that look like any of the $T_{U_i}$ before merging, in such a manner that after merging they collectively reproduce $T_{U_1 \cdots U_n}$. This prompts us to define a class $[T_{U_1} \cdots T_{U_n}]$ of all those towers $T'_{U_1}, \cdots, T'_{U_n}$ such that $\oT_{1 \leq i \leq n} T'_{U_i} = T_{U_1 \cdots U_n}$.

\subsection{Category $\TqBunp$}
Having introduced $\cT_q$ above, we can consider the functor:
\begin{align}
	\chi_{p,q}: \bOne[\peps]^q & \rarr \Grpd \nonumber \\
	(U_1,\cdots, U_q) & \mapsto  \cT_{U_1 \cdots U_q} = \{ T_{U_1 \cdots U_q} \} \nonumber
\end{align}
with $ T_{U_1 \cdots U_q} \doteq T_{U_1} \oT \cdots \oT T_{U_q}$. Given $f = f_1 \times \cdots \times f_q: (U_1, \cdots, U_q) \rarr (U'_1, \dots, U'_q)$ in the base, the functor $\chi_{p,q}$ produces a map $\chi_{p,q}f$ of groupoids, properly defined as $\chi f_1 \oT \cdots \oT \chi f_q$ as a functor, as a result of each $\chi f_i$ being a functor, and $\oT^q$ being a multi-functor.\\

Doing a Grothendieck construction on $\chi_{p,q}$ yields a bundle $\TqBunp$ over $\bOne[\peps]^q$, with objects of the form $((U_1, \cdots, U_q),T_{U_1 \cdots U_q})$, and a morphism of such an object to another $((U'_1, \cdots, U'_q), T_{U'_1 \cdots U'_q})$ is given by $(f,\alpha)$, $f$ defined above, $\alpha = \oT_{1 \leq i \leq q} \alpha_i$, $\alpha_i: \chi_p f_i \xi \cong \psi$, $1 \leq i \leq q$, for $\xi \in \cT_{U_i}$, $\psi \in \cT_{U'_i}$. The forgetful functor shows that $\TqBunp$ is cofibered in groupoids over $\bOne[\peps]^q$.

\subsection{$\bK$-formalism}
\subsubsection{$\delta$ map on towers}
Collecting all those bundles defined above together, we form the categorical bundle $\TdBund$ cofibered in groupoids grade-wise that is $\TdBund = \{\TqBunp\}_{p,q}$ with $\TqBunp  \rarr \bOne[\peps]^q$ cofibered in groupoids. Recall that in the $\delta$-formalism we had the delta maps $\cdots \bOne[\peps] \xrarr{\delta} \bOne[\pponeeps] \rarr \cdots$ acting on thick points. Above points we have towers, which are nothing but stacked $\bK$-neighborhoods, themselves made of points. If one can thicken points, one can thicken collections of points, hence cross-sections of towers, hence the towers themselves, something we used implicitly in constructing towers anyway. Thus over thickened elements $\delta U$ for $U \in \bOne[\peps]$, one can have thickened towers $\delta T_U$ over $\bOne[\pponeeps]$ that still contain $T_U$ as their core. We can see this as having a vertical $\bK$-dimension, and a horizontal $\mu R$-dimension to account for thickening. Hence we get a map $\delta: \TBunp \rarr \TBun_{p+1}$, which we will show is functorial, and which we can further generalize to a functor $\delta_q: \TqBunp \rarr \TqBun_{p+1}$ equivariant with respect to the tensor product $\oT$ on $\cT_{\centerdot}$. Letting $\delta = \oP_{q \geq 1} \delta_q$, $(\TdBund, \oT, \delta)$ defines what we call the $\bK$-formalism.\\

As an aside, each tower ends in $\Omega$, which is invariant, thus operating $\delta$ on each tower thickens it, except $\Omega$, its limit, which remains invariant. One can see that this is indeed the case using the following argument. $\Omega$ being part of the system $T_U$ (it being its limit), it follows $\delta \Omega$ is part of $ \delta T_U$, hence there exists a map $\Omega = \lim \delta T_U \rarr \delta \Omega$. If we have no equality, there exists $V \in T_U$ such that $\delta V = \delta \Omega$, but it would then follow that $V \rarr \Omega$ is beyond $\lim T_U$, which is not possible.\\

\subsubsection{categorical formalism}

-\underline{$\delta$ functors}\\
We show $\delta$ is a functor on $\TBun$. Consider $\TBunp$ over $\bOne[\peps]$, $(U,T_U)$ an object thereof. We have $\delta (U,T_U) = (\delta U, \delta T_U)$, $\delta T_U$ defined above. For a morphism $(f,\alpha): (U,T_U) \rarr (V,T_V)$, we have an associated map $\delta f: \delta U \rarr \delta V$. Regarding the thickened towers, if we have such thickened maps on thick points, we can have the same for thickened open sets, hence for cross-sections of towers, hence for towers themselves. Thus if we have a map $\alpha: \chi (f) T_U \rarr T_V $ of towers, we can construct its extension $\delta \alpha: \delta (\chi (f) T_U) \rarr \delta T_V$. More specifically, recall how $\alpha$ is defined:
\beq
\xymatrix{
	T_U \ar@{.>}[dr] \ar[r]^-{\chi f} & \chi f (T_U) = T'_V \ar[d]^{\alpha} \\
&T_V
} \nonumber
\eeq
When we write $\alpha$, we really mean the composite $\alpha \circ \chi f$. All those maps above are maps of towers. We can thicken such maps, hence $\delta \alpha$ is well-defined. Finally $\delta(f,\alpha) = (\delta f, \delta \alpha)$. Because $\delta$ is a functor, $\delta(id,id) = (id,id)$, and we also have composition: 
\begin{align}
	\delta((g,\beta) \circ (f,\alpha))&= \delta ((g \circ f, \beta \circ \alpha)) \nonumber \\
	&= (\delta( g \circ f), \delta( \beta \circ \alpha)) \nonumber \\
	&= (\delta g \circ \delta f, \delta \beta \circ \delta \alpha) \nonumber \\
	&= (\delta g, \delta \beta) \circ (\delta f, \delta \alpha) \nonumber \\
	&= \delta (g, \beta) \circ \delta (f,\alpha) \nonumber
\end{align}
by virtue of the fact that $\delta(\beta \circ \alpha) = \delta \beta \circ \delta \alpha$, which follows from working cross-section-wise on towers and using the fact that $\delta$ is functorial on neighborhoods. Here $\beta \circ \alpha$ is short for:
\beq
\xymatrix{
	T_U \ar@{.>}[dr] \ar[r]^-{\chi f} &\chi f T_U = T'_V \ar[d]^{\alpha} \\
	&T_V \ar@{.>}[dr] \ar[r]^-{\chi g} & \chi g T_V = T'_W \ar[d]^{\beta} \\
	&& T_W
} \nonumber
\eeq

-\underline{$\delta_q$ functors}\\
We now show:
\beq
\xymatrix{
	\delta_q: \TqBunp \ar[d] \ar@{.>}[r] & \TqBun_{p+1} \ar[d] \\
	\dOne[\peps]^q \ar[r]_{\delta^q} & \dOne[\pponeeps]^q
} \nonumber
\eeq
is a functor as mentioned above. Now that we have towers merging, one has to be careful about thickening merged towers. For a generic multi-tower $\TUoneUq$ suppose one can find towers $T_{U_1}, \cdots, T_{U_q}$ such that one can write $\TUoneUq = \oT_{1 \leq i \leq q} T_{U_i}$. We want the thickening of $\TUoneUq$ to be independent of this representation. Let $T'_{U_1} \cdots T'_{U_q}$ be towers in the class $[T_{U_1} \cdots T_{U_q}]$. Observe that $i$ being fixed, if $T'_{U_i}$ and $T_{U_i}$ look alike before merging, it follows their thickenings are the same. After merging, points of $T'_{U_1} \oT \cdots \oT T'_{U_q}$ and $T_{U_1} \oT \cdots \oT T_{U_q}$ coincide so thickening such points is independent of the representation. Thus overall we have  $ \delta \TUoneUq = \oT \delta T_{U_i} = \oT \delta T'_{U_i}$ is independent of the choice of representative upon thickening. Thus $\delta_q$ is well defined on $\TqBunp$ on objects.\\

Now if $(U_1, \cdots, U_q)$ and $(V_1, \cdots, V_q)$ are object in the base with $f_i:U_i \rarr V_i$ for $1 \leq i \leq q$, one has $\delta^q(U_1, \cdots, U_q) = (\delta U_1, \cdots, \delta U_q)$, and $\delta^q (f_1, \cdots, f_q) = (\delta f_1, \cdots, \delta f_q)$. Then $\delta^q$ is a functor by virtue of the fact that $\delta$ itself is a functor, and $\times^q$ is a multi-functor. The same is true of $\oT^q$, which makes $\delta_q := \oT^q \delta$ with $\delta$ on bundles a functor. Finally letting $\delta_{\cdot} = \oP_{q \geq 1} \delta_q$, $\delta: \Td \Bun_p \rarr \Td \Bun_{p+1}$ is a functor grade-wise in $q$, which makes $(\TdBund, \oT, \delta)$ a categorical directed system that we define as the \textbf{$\bK$-formalism}.

\section{Duality of formalisms}
\subsection{Active/passive formalisms}
Observe that the $\delta \Phi$ provides us with a quantized, horizontal deformation of neighborhoods/thick points, by construction. Writing $\delta \Phi \doteq \dOne =  \Omega$ this can also be viewed as providing a progressive immersion of the physical space $\dOne$ into the ambient vista $\Omega$. In contrast, the $\bK$-formalism (abbr. $\bK \Phi$) provides us with a smooth, vertical deformation of neighborhoods, presented by thickened towers over the thick points of the natural $\delta \Phi$. Essentially the $\delta \Phi$ corresponds to a passive, spontaneous immersion of $\dOne$ within $\Omega$, whereas the $\bK \Phi$ presents an active $\bK$-enlargement of $\bOne$, active insofar as it is $\bK$-generated, resulting from a prismatic decomposition of $\delta \Phi$-neighborhoods into their $\bK \Phi$ components from the perspective of the $\delta \Phi$.

\subsection{Global picture}
Note that we heuristically have:
\beq
r \Omega \doteq rec \delta \Phi \leadsto \bK \Phi \nonumber
\eeq
By this we mean the inert space $\Omega$ is dynamically presented in the guise of $\delta \Phi$-resolutions which one can regard as flows. We remark in passing that in a single thick point $U \in \dOne[\peps]$, one has the totality of $\Omega$, the reason being that each standalone $\dOne[\peps]$ covers $\Omega$, and each object $U$ of $\dOne[\peps] \in \Cat$ is linked to every other object therein. Those $\delta \Phi$-resolutions are then functorially mapped to the towers used in the $\bK \Phi$. That last transition $\delta \Phi \rarr \bK \Phi$ we regard as a prismatic decomposition, which one would assume entails a loss of information since we decompose pure knowledge in $\bOne = r \Omega$ into consciousness plus obscurations, and the transition $\delta \Phi \leadsto \bK \Phi$ focuses on $\bK$ alone. The reverse move $\bK \Phi \rarr \delta \Phi$ is produced via the Grothendieck construction. Indeed we have seen that functors $\chi_{p,q}: \bOne[\peps]^q \rarr \Grpd$ collectively define categories $\TqBunp \rarr \bOne[\peps]^q$ cofibered in groupoids (via the Grothendieck construction) thus we have $\bK \Phi \leadsto \delta \Phi$. A natural question is whether we have a relation between $\bK \Phi$ and $\delta \Phi$. We will show that we have $r\Omega \doteq \bK \Phi$, thereby showing that there is actually no loss of information, and that we have a duality $\delta \Phi \sim \bK \Phi$.\\

\begin{Kay}
	One thing we have to specify at this point is what exactly does ``consciousness" mean in our context. So far we have identified this with awareness in the pedestrian sense, which naively results from existence. But not everything that exists is endowed with awareness. Existence is a preeminent concept, the only characteristic if any of $\Omega$. We equate (pure) $\bK$ with existence, which can be treated categorically. Thus what we call (impure) $\bK$ at our level is only a mind-related type of interaction, which we just refer to as awareness. In this manner consciousness appears as a truncation of pure $\bK$. 
\end{Kay}

\subsection{Presentation of dual pictures}
$\Omega$ is inert, but can be dynamically presented by its covers provided by the $\delta \Phi$. In a first time, $\delta \Phi \doteq \Omega$. \\

On the other hand we also have the $\bK \Phi$. To have a duality $\bK \Phi \sim \delta \Phi$ would imply the Grothendieck construction is an equivalence in this particular situation.\\

\subsection{Sums of $\bK$}
Now given $f:U \rarr V$ in the base $\bOne[\peps]$, $\chi_p f:\cT_U \rarr \cT_V$ is a functor that maps individual towers to towers. Given a path $\gamma$ in the base $\bOne[\peps]$, we therefore have paths in towers which we call threads, and which we can collectively denote by $\bK(\gamma)$, those are all the possible threads that can be linked to $\gamma$ in the base. Using the Grothendieck construction, we have a forgetful functor: $\pi: \TBunp \rarr \bOne[\peps]$ for which morphisms are given between random towers, not just those that are along a same path, a same thread. Thus this heuristically corresponds to $\sum \bK(a)$ in full.\\

\subsection{$r\Omega$ as sum of $\bK$}
Next we argue $r\Omega  = \sum \bK(a)$. Again we certainly have an inclusion to the left. We claim this is actually an equality. The sum on the right is over $a \in \bOne \doteq r\Omega$. Even if it is true that over points $a \in \bOne$, towers give us incomplete $\bK$, the problem being that none of the cross-sections of any tower can cover any region of $\Omega$, all towers end in $\Omega$ by definition thus eventually $\sum \bK(a)$ provides a cover of $\Omega$ $\bK$-wise. Indeed that all towers have $\Omega$ for limit was by design, and this is not hard to see; even if $\tilde{U}$ is large enough to cover $r \Omega$, one can still consider higher cross-sections whose local intensities are higher than those of $\tilde{U}$. Thus ultimately in the limit we consider pure $\bK$, which is none other than $r \Omega$. It is worth emphasizing this change of perspective, from a point-base formalism provided by finite towers, to one where referents are absent, namely $\Omega$, a transcendental transition akin to a ``cloud of unknowing". Dynamically, high up in the towers, $\iota(a) \ll \iota(\Omega)$, at which point things are no longer centered, concepts dissolve, and obscurated consciousness turns into pure consciousness.

\subsection{Duality}
To conclude, in the $\delta \Phi$, $\mu R$ dictates the dynamics, it is a globally $\Omega$-induced phenomenon, measured locally, and we eventually get $\dOne \doteq  \Omega$ (or equivalently $\bOne \doteq r \Omega$). On the $\bK$-side/reconstructed side of things, we deal with the natural realm, with obscurations, hence local consciousness now is relevant since global pure $\Omega$-level consciousness is lost. Working locally we eventually arrive at $\bK \Phi \doteq r\Omega$. Thus $\delta \Phi$ is globally induced, $\bK \Phi$ is locally generated, and both recover $\Omega$. This forms a heuristic T-duality:
\beq
\xymatrix{
	&\Omega \ar[dl]_{\text{flow}} \ar[dr]^{rec}_{\cong} \\
	\Omega \doteq \delta \Phi && r \Omega \doteq \bK \Phi
} \nonumber
\eeq

\begin{inc}
	Precisely because $a$ never reaches $\Omega$ via any tower $T_a$, albeit in the limit, this is just a formal duality. As a consequence of this incompleteness however, it follows the practical, real life $\bK \Phi \varsubsetneq r \Omega$, or in other terms $\sum \bK(a)$ is not sufficient to reproduce $\Omega$ in full.
\end{inc}

\section{Determinism vs free will}
In the present work we have argued $\dOne$ can be immersed into $\Omega$, the idea being that one can categorically identify which events or close to one another, which subsumes an underlying deterministic picture of $\dOne$, or equivalently $\bOne$, where natural phenomena take place. This is also clearly evident if for $\bOne$ one adopts a representation by the Segal category $\delta^{\infty}\cX$. On the other hand in the $\bK \Phi$, we discuss consciousness, which has a strong association with the concept of free will. To be fair, the exact connection remains elusive and is far from being settled. For our purposes, we take free will to designate spontaneous generations from the perspective of the agent free will is in reference to, broadly speaking. No free will implies no possible consciousness since we would just have a flow of information through the agent. By contrapositive, to have consciousness subsumes there is such a thing as free will. If it was, and still is, difficult to define consciousness, not to have free will as underlying assumption is an even more serious problem, notwithstanding the fact that it seems $\bK$ is defined within a deterministic ambient environment $\bOne$.\\

We argue there is actually no contradiction in having this picture in mind, and that we can actually have both concepts coexisting. The argument goes as follows. Consider a system endowed with free will, that we can represent by a category $\cC$. Its manifestations in the realm of natural phenomena can be represented by a functor $F: \cC \rarr \Catinf$. In particular free will is a coherent manifestation in $\bOne$, so must be represented by such a functor. However by the $\infty$-categorical Grothendieck construction (\cite{Lu}), discussed at length in \cite{RG6}, we have $\Fun(\cC,\Catinf) \simeq \coCartFib(\cC)$. Thus to any functor $F: \cC \rarr \Catinf$ corresponds a co-cartesian fibration over $\cC$. This is a functor into $\cC$, hence is purely deterministic. In this duality, $\cC$ is a pivot, and from that perspective the LHS of this duality presents $\Fun(\cC,\Catinf)$ as providing something with no precedent prior to $\cC$, it being the point of focus, thus we deal with something that can be interpreted as being free will-based. As argued in \cite{RG6}, this duality provides an algebraic flow of information through $\cC$, we have a conservation law in a sense, hence free will can be seen as being deterministically produced. Additionally, what this duality says is that it is perpetual, and simultaneous; information flowing into $\cC$ generates functors that correspond to manifestations of $\cC$ in its ambient universe, which in turn are fed back into co-cartesian fibrations over $\cC$, ad infinitum.\\

The individual decision making process also results from dynamics within $\cC$, and this also is encoded in this Grothendieck duality; external events shape $\cC$ via co-cartesian fibrations, which generate functors $F: \cC \rarr \Catinf$, some of which produce internal $\infty$-categories, read self-characteristics, some of which will be relevant for generating awareness. So far we are saying our inner self is produced in a deterministic fashion. Those internal functors are then fed back into co-cartesian fibrations to produce other functors, some of which will produce a basis for free will. Thus our own self can be argued to be formed in a deterministic fashion, and those newly formed characteristics can contribute to our decision making process. Observe also that the inner mechanics of it all is deterministic but unnoticeable since it arises from myriads of fibrations $\cD \rarr \cC$ that generate inner functors $F: \cC \rarr \Catinf$ that shape our inner dynamics, which ultimately leads in particular to free will at a global level at which scale the inner deterministic mechanics is not apparent.\\

To make this duality determinism-free will more clear, recall from \cite{RG6} that both theories can be deconstructed into a common ground state, in which both theories are in deconstructed form, and one can either reconstruct the ground state into co-cartesian fibrations, which is deterministic, or into functors from $\cC$ into $\Catinf$, which subsumes $\cC$ has free will. However in the ground state itself those concepts are in dissolved form and therefore become moot. In other terms ultimately it does not matter whether things are deterministic or spontaneous. Another argument in favor of the non-importance of this apparent issue, is that it originates from identifying a fixed category $\cC$, which is a point of reference, and introducing referents necessarily entails dualities, and therefore possible incompatibilities, case in point the determinism/free will conundrum. Those are relative concepts only.

\bigskip
\footnotesize
\noindent
\textit{e-mail address}: \texttt{rg.mathematics@gmail.com}.

\end{document}